\documentclass[12pt]{article}
\usepackage{amsmath}
\usepackage{amsfonts}
\usepackage{mitpress}
\usepackage{verbatim}
\usepackage{graphicx}
\usepackage{caption}
\usepackage{subcaption}
\usepackage[bookmarks=true,         bookmarksnumbered=true,
colorlinks=true,
pdfstartview=FitV,
linkcolor=blue, citecolor=blue, urlcolor=blue,
]{hyperref}
\usepackage[sort,compress]{cite}

\usepackage{soul}
\usepackage{array}

\newcommand{\kpar}[1]{\bigskip\noindent{\em {#1}}}
\newcommand{\floor}[1]{{\lfloor{#1}\rfloor}}

\begin{document}

\date{Last round of revision: \\
September 24, by AC and FL\\
September 21, by KL \\
September 17, by AC and FL \\
September 13, by FL\\
September 9, by KL\\
August 23, by AC}


\title{Data-based stochastic model reduction for the
  Kuramoto--Sivashinsky equation}

\author{Fei Lu\thanks{Department of Mathematics, University
    of California at Berkeley and Lawrence Berkeley National
    Laboratory. E-mail addresses: feilu@berkeley.edu (F. Lu, corresponding author); chorin@math.berkeley.edu (A.J. Chorin)}, Kevin K.~Lin\thanks{School of Mathematics, University
    of Arizona. E-mail address: klin@math.arizona.edu}, and Alexandre J.~Chorin$^*$}
\maketitle

\begin{abstract}

The problem of constructing data-based, predictive, reduced models for the Kuramoto-Sivashinsky equation is considered, under circumstances where one has observation data only for a small subset of the dynamical variables.  Accurate prediction is achieved by developing a discrete-time stochastic reduced system, based on a NARMAX (Nonlinear Autoregressive Moving Average with eXogenous input) representation. The practical issue, with the NARMAX representation as with any other, is to identify an efficient structure, i.e., one with a small number of terms and coefficients.  This is accomplished here by estimating coefficients for an approximate inertial form.  The broader significance of the results is discussed. 
 \end{abstract}


\section{Introduction}

There are many high-dimensional dynamical systems in science and
engineering that are too complex or computationally expensive to solve
in full, and where only a relatively small subset of the degrees of
freedom are observable and of direct interest.
%
%
Under these conditions, it is useful to derive low-dimensional models
that can predict the evolution of the variables of interest without
reference to the remaining degrees of freedom, and reproduce their
statistics at an acceptable cost.

We assume here that the variables of interest have been observed in the
past, and we consider the problem of deriving low-dimensional models on
the basis of such prior observations. We do the analysis in the case of
the Kuramoto--Sivashinsky equation (KSE):
\begin{eqnarray}
&&\frac{\partial v}{\partial t}+v\frac{\partial v}{\partial x}+\frac{\partial ^{2}v}{\partial x^{2}}+\frac{\partial ^{4}v}{\partial x^{4}}=0\,,x\in \mathbb{R},t>0;  \label{KSE} \\[1ex]
&&v(x,t)=v(x+L,t);\,\,\,v(x,0)=v_{0}(x),  \notag
\end{eqnarray}
where $t$ is time, $x$ is space, $v$ is the solution of the equation,
$L$ is an assumed spatial period, and $v_0$ is the initial datum.  We
pick a small integer $K$, and assume that one can observe only the
Fourier modes of the solution with wave numbers $k=1,\dots K$ at a
discrete sequence of points in time.  To model the usual situation where
the observed modes are not sufficient to determine a solution of the
differential equations without additional input, we pick $K$ small
enough so that the dynamics of a Galerkin-Fourier representation of the
solution, truncated so that it contains only $K$ modes, are far from the
dynamics of the full system.  The goal is to account for the effects of
``model error'', i.e. for the effects of the missing ``unresolved" modes
on the ``resolved" modes, by suitable terms in reduced equations for the
resolved modes, using the information contained in the observations of
the resolved modes; we are not interested in the unresolved modes per se. In the
present paper the observations are obtained by a solution of the full
system; we hope that our methods are applicable to problems where the
data come from physical measurements, including problems where a full
model is not known.

We start from the truncated equations for the resolved modes, and solve
an inverse problem where the data are used to estimate the effects of
model error, i.e., what needs to be added to the truncated equations for
the solution of the truncated equations to agree with the data.  Once
these effects are estimated, they need to be identified, i.e.,
summarized by expressions that can be readily used in computation. In
our problem the added terms can take a range of values for each value of
the resolved variables, and therefore a stochastic model is a better
choice than a deterministic model. We solve the inverse problem within a
discrete-time setting (see \cite{CL15}), and then identify the needed
terms within a NARMAX (Nonlinear Autoregression with Moving Average and
eXogenous input) representation of discrete time series.  The
determination of missing terms from data is often called a
``parametrization"; what we are presenting is a discrete stochastic
parametrization.  The main difficulty in stochastic parametrization, as
in non-parametric statistical inference problems, is making the
identified representation efficient, i.e., with a small number of
  terms and coefficients. We accomplish this by a
semi-parametric approach: we propose terms for the NARMAX representation
from an approximate ``inertial form'' \cite{JKT90}, i.e., a system
  of ordinary differential equations that describes the motion of the
  system on a finite-dimensional, globally attracting manifold called an
  ``inertial manifold". (Relevant facts from
inertial manifold theory are reviewed later.)

A number of stochastic parametrization methods have been proposed in
recent years, often in the context of weather and climate prediction. In
\cite{Pal01, Wil05,AMP13}, model error is represented as the sum of an
approximating polynomial in the resolved variables, obtained by
regression, and a one-step autoregression. The shortcomings of this
representation as a general tool are that it does not allow the model
error to depend sufficiently on the past values of the solution, that
the model error is calculated inaccurately, especially when the data are
sparse, and that the autoregression term is not necessarily small,
making it difficult to solve the resulting stochastic equations
accurately. Detailed comparisons between this approach and a
discrete-time NARMAX approach can be found in \cite{CL15}.  In
\cite{CVE08,Kwa12} the model error is represented as a conditional Markov
chain that depends on both current and past values of the solution; the
Markov chain is deduced from data by binning and counting, assuming that
exact observations of the model error are available, i.e., that the
inverse problem has been solved perfectly. It should be noted that the
Markov chain representation is intrinsically discrete, making this work
close to ours in spirit. In \cite{MH13} the noise is treated as
continuous and represented by a hypo-elliptic system that is partly analogous to the NARMAX representation, once translated from the continuum to the
grid. An earlier construction of a reduced approximation can be found in
\cite{CH14}, where the approach was not yet fully discrete.  Other
interesting related work can be found in \cite{CKG11, DN07,DD09, Har16, KCG15}.
The present authors' previous work on the discrete-time approach
to stochastic parametrization and the use of NARMAX representations can
be found in \cite{CL15}.

The KSE is a prototypical model of spatiotemporal chaos.  As a
  nonlinear PDE, it has features found in more complex models of
  continuum mechanics, yet its analysis and numerical solution are
  fairly well understood because of its relatively simple
  structure. There is a lot of previous work
on stochastic model reduction for the KSE.  Yakhot\cite{Yak81} developed
a dynamic renormalization group method for reducing the KSE, and showed
that the model error generates a random force and a positive
viscosity.  Recent development of this method can be found in \cite{SPKP13, SPPK15}. Toh \cite{Toh87} studied the statistical properties of the
KSE, and constructed a statistical model to reproduce the energy
spectrum. Rost and Krug \cite{RK95} presented a model of interacting
particles on the line which exhibits spatiotemporal chaos, and made a
connection with the stochastic Burgers equation and the KPZ equation.

Stinis \cite{Sti04} addressed the problem of reducing the KSE as an
under-resolved computation problem with missing initial data, and used
the Mori-Zwanzig (MZ) formalism \cite{Zwan01} in the finite memory approximation
\cite {CHK02, CH13} to produce a reduced system that can
make short-time predictions. Full-system solutions were used to compute
the conditional means used in the reduced system. As discussed in
\cite{CL15}, the NARMAX representation can be thought of as both a
generalization and an implementation of the MZ formalism, and the
full-system solutions used by Stinis can be viewed as data, so that the
pioneering work of Stinis is close in spirit to our work. We provide
below a comparison of our work with that of Stinis.

\bigskip

The paper is organized as follows. In section 2, we introduce the
Kuramoto--Sivashinsky equation, the dynamics of its solutions, and its
numerical solution by spectral methods. In section 3 we apply the
discrete approach for the determination of reduced systems to the KSE,
and discuss the NARMAX representation of time series. In section 4 we
use an inertial form to determine the structure of a NARMAX representation for the KSE, and estimate its coefficients.  Numerical results are presented in section 5. Conclusions
and the broader significance of the work, as well as its limitations,
are discussed in a concluding section.

\section{The Kuramoto-Sivashinsky equation}
\label{sect:kse}

We begin with basic observations: in Eq.~(\ref{KSE}), the term
${\partial^{2}v}/{\partial x^{2}}$ is responsible for instability at
large scales, the dissipative term ${\partial ^{4}v}/{\partial x^{4}}$
provides damping at small scales, and the non-linear term $v{\partial
  v}/{\partial x}$ stabilizes the system by transferring energy between
large and small scales.  To see this, first write the KSE in terms of
Fourier modes: 
\begin{equation}
\frac{d}{dt}v_{k}=(q_{k}^{2}-q_{k}^{4})v_{k}-\frac{iq_{k}}{2}\sum_{l=-\infty
}^{\infty }v_{l}v_{k-l},  \label{FM}
\end{equation}
where the $v_{k}(t)$ are the Fourier coefficients
\begin{equation}
v_{k}(t)~:=~\mathcal{F}[v(\cdot
  ,t)]_{k}~:=~\frac{1}{L}\int_{0}^{L}v(x,t)e^{-iq_{k}x}dx,
\end{equation}
where $q_{k}=\frac{2\pi k}{L}$, $k\in{\mathbb Z}$, and ${\mathcal F}$
denotes the Fourier transform, so that
\begin{equation}
v(x,t)=\mathcal{F}^{-1}[v_{\cdot }(t)]=\sum_{k=-\infty }^{+\infty }v_{k}(t)e^{iq_{k}x}.
\end{equation}
Since $v$ is real, the Fourier modes satisfy $v_{-k}=v_{k}^{\ast }$,
where $v_{k}^{\ast }$ is the complex conjugate of $v_{k}$.  We refer to
$|v_k(t)|^2$ as the ``energy" of the $k$th mode at time $t$.

\begin{figure}
\centering
\includegraphics{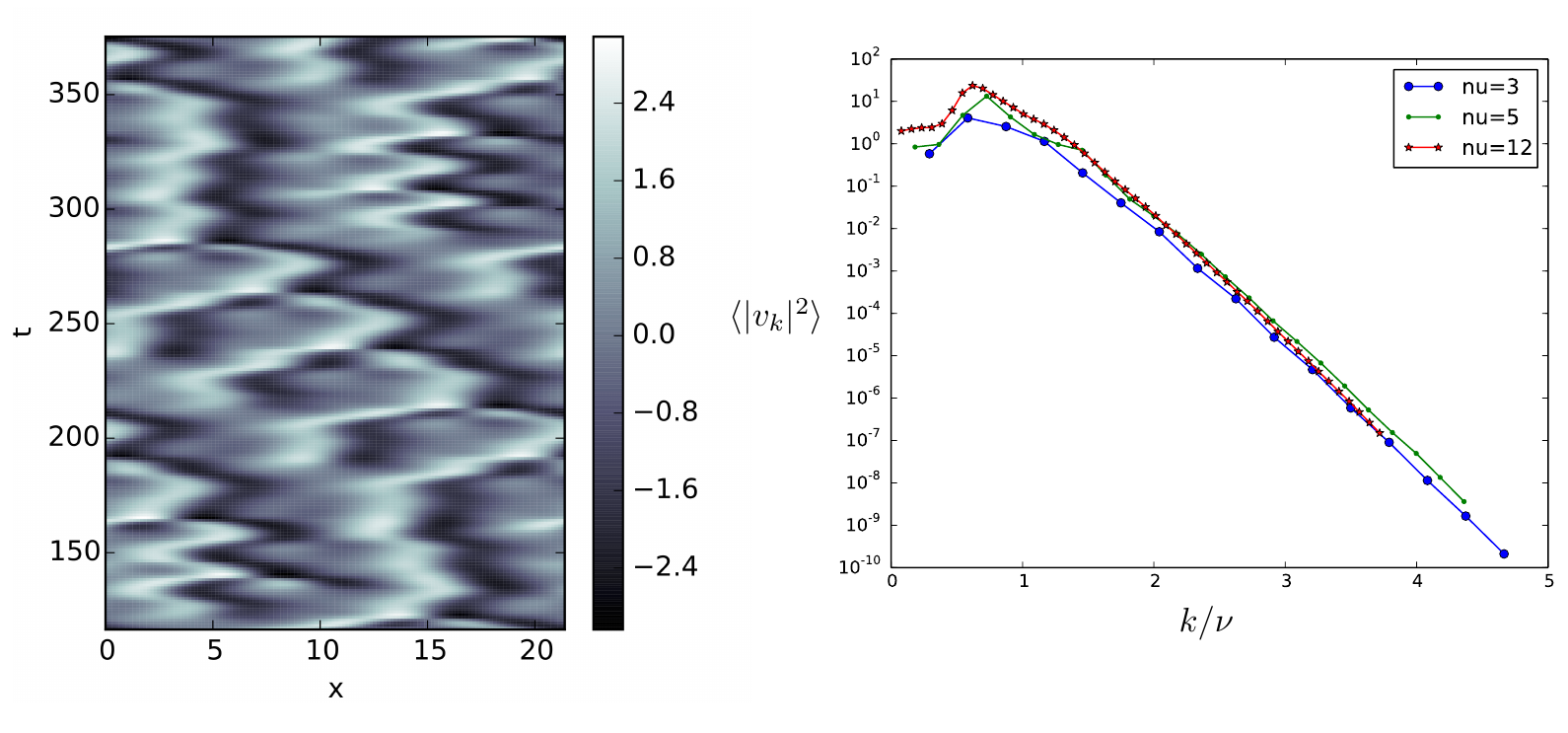}
  \caption{Solutions of the KSE.  {\em Left:} a sample solution of the
    KSE.  {\em Right:} Mean energy $\langle|v_k|^2\rangle$ as a function
    of $k/\nu$, where $\nu=L/2\pi$ is the number of unstable modes, for
    a number of different values of $L$.  }
  \label{fig:kse-solutions}
\end{figure}

Next, we consider the linearization of the KSE about the zero solution. In the linearized
equations, the Fourier modes are uncoupled, each represented by a
first-order scalar ODE with eigenvalue $q_k^2-q_k^4$.  Modes with
  $|q_{k}|>1$ are linearly stable; modes with $|q_k|\leq1$ are not.  The
  linearly unstable modes, of which there are $\nu=\floor{L/2\pi}$, are
coupled to each other and to the damped modes through the nonlinear
term.  Observe that if the nonlinear terms were not present, most
initial conditions would lead to solutions whose energies blow up
exponentially in time.  The KSE is, however, well-posed (see, e.g.,
\cite{CFNT88}), and it can be shown that solutions remain globally
bounded in time (in suitable function spaces) \cite{Goo94, BG06}.  The solutions of
Eq.~(\ref{KSE}) do not grow exponentially because the quadratic
nonlinearities, which formally conserve the $L^2$ norm, serve to
transport energy from low to high modes.  Figure~\ref{fig:kse-solutions}
shows an example solution, as well as the {\em energy spectrum}
$\langle|v_k|^2\rangle,$ where $\langle \phi\rangle$ denotes the limit of
$\frac1T\int_0^T \phi(v(t))dt$ as $T\to\infty~.$ The energy spectrum has
a characteristic ``plateau'' for small wave numbers, which gives way to
rapid, exponential decay as $k$ increases, see Figure~\ref{fig:kse-solutions}(right). This concentration of   energy in the low-wavenumber modes is reflected in the cellular
  character of the solution in Figure~\ref{fig:kse-solutions}(left),
  where the length scale of the cells is determined by the modes 
  carrying the most energy.
  
Another feature of the nonlinear energy transfer is that the KSE
possesses an inertial manifold $M$.  When the KSE is viewed
as an infinite-dimensional dynamical system on a suitable function space
$B$, there exists a finite-dimensional submanifold $M\subset B$ such
that all solutions of the KSE tend asymptotically to $M$ (see
e.g. \cite{CFNT88}).  The long-time dynamics of the KSE are thus
essentially finite-dimensional.  Moreover, it can be shown that for
sufficiently large $N$ ($N>\mbox{dim}(M)$ at the very minimum),
an inertial manifold $M$ can be written in the form
\begin{equation}
  M = \{x +\psi(x)|x\in\pi_N(B)\}
\end{equation}
where $\pi_N$ is the projection onto the span of $\{e^{ikx},
k=1,\cdots,N\}$, and $\psi:\pi_N(B)\to(\pi_N(B))^\perp$ is a
Lipschitz-continuous map.  That is to say, for trajectories lying on the
inertial manifold $M$, the high-wavenumber modes are completely
determined by the low-wavenumber modes.  Inertial manifolds will be
useful in Section~\ref{sect:narmax-struct}, and we say more about them
there.

The KSE system is Galilean invariant; if $v(x,t)$ is a solution, then
$v(x-ct,t)+c$, with $c$ an arbitrary constant velocity, is also a
solution.  Without loss of generality, we set $\int v(x,0)dx=0$, which
implies that $v_{0}(0)=0$.  From (\ref{FM}), we see that $v_{0}(t)\equiv
0$ for all $t$ and $\int v(x,t)dx\equiv 0$.  In physical terms,
solutions $v(x,t)$ of the KSE can be interpreted as the velocity of a
propagating ``front,'' for example as in flame propagation, and this
decoupling of the $v_0$ equation from $v_k$ for $k\neq0$ means the mean
velocity is conserved.

\kpar{Chaotic dynamics and statistical assumptions.} Numerous studies
have shown that the KSE exhibits chaotic dynamics, as characterized by a
positive Lyapunov exponent, exponentially decaying time correlations, and other signatures of chaos (see, e.g., \cite{HN86}
and references therein).  Roughly speaking, this means that nearby
trajectories tend to separate exponentially fast in time, and that,
though the KSE is a deterministic equation, its solutions are
unpredictable in the long run as small errors in initial conditions are
amplified.  Chaos also means that a statistical modeling approach is
natural.  In what follows, we assume our system is in a chaotic regime,
characterized by a translation-invariant physical invariant probability
measure (see, e.g., \cite{ER85,Y13} for the notion of physical invariant
measures and their connections to chaotic dynamics).  That is, we assume
that numerical solutions of the KSE, sampled at regular space and time
intervals, form (modulo transients) multidimensional time series that
are stationary in time and homogeneous in space, and that the resulting
statistics are insensitive to the exact choice of initial conditions.
Except where noted, this assumption is consistent with numerical
observations.  Hereafter we will refer to this as ``the'' ergodicity
assumption, as this is the assumption of ergodicity needed in the
present paper.

The translation-invariance part of our ergodicity assumption has
  the consequence, the Fourier coefficients cannot have a preferred
  phase in steady state, i.e., the physical invariant measure has the
  property that if we write the $k$th Fourier coefficient $v_k$ (see
  below) as $a_ke^{i\theta_k}$, then the $\theta_k$ are
  uniformly distributed.  The theory of stationary stochastic processes
  also tells us that distinct Fourier modes are uncorrelated
  \cite{CH13}, though they are generally not independent: one can
readily show this by, e.g., checking that ${\rm
  cov}(|v_k|^2,|v_\ell|^2)\neq0$ numerically.  (Such energy correlations
are yet another consequence of the nonlinear term acting to transfer
energy between modes.)


\kpar{Numerical solution.}  For purposes of numerical approximation, the
system can be truncated as follows: the function $v(x,t)$ is sampled at
the grid points $x_{n}=nL/N, n=0,\dots ,N-1$, so that $v^{N}=\left(
v(x_{0},t),\dots ,v(x_{N-1},t)\right) $, where the superscript $N$ in
$v^{N}$ is a reminder that the sampled $v$ varies as $N$ changes. The
Fourier transform $\mathcal{F}$ is replaced by the discrete Fourier
transform ${\mathcal F}_{N}$ (assuming $N$ is even):
\begin{equation*}
v_{k}^{N}={\mathcal F}_{N}[v(\cdot
,t)]_{k}=\sum_{n=0}^{N-1}v^{N}(x_{n},t)e^{-iq_{k}x_{n}},\,
\end{equation*}
and
\begin{equation*}
\,v^{N}(x_{n},t)={\mathcal F}_{N}^{-1}[v_{\cdot }^{N}(t)]=\frac{1}{N}
\sum_{k=-N/2+1}^{N/2}v_{k}^{N}(t)e^{iq_{k}x_{n}}.
\end{equation*}
Noting that $\hat{v}_{0}^{N}=0$ due to Galilean invariance, and
setting $\hat{v}_{N/2}^{N}=0$, we obtain a truncated system
\begin{equation}
\frac{d}{dt}v_{k}^{N}=(q_{k}^{2}-q_{k}^{4})v_{k}^{N}-\frac{iq_{k}}{2}\sum
_{\substack{ 1\leq |l|\leq N,  \\ 1\leq \left\vert k-l\right\vert \leq N}}
v_{l}^{N}v_{k-l}^{N},~  \label{DFM}
\end{equation}
with $k=-N/2+1,\dots ,N/2$. This is a system of ODE (of the
real and imaginary parts of $v_{k}^{N}$) in $\mathbb{R}^{N-2}$, since
$v_{-k}^{N}=v_{k}^{N,\ast }$ and $v_{N/2}^{N}=$ $v_{0}^{N}\equiv 0$.

Fourier modes with large wave numbers are typically small and can be
neglected, as can be seen from a linear analysis: $v_{k}$ decreases at
approximately the rate $e^{(q_{k}^{2}-q_{k}^{4})t}$, where $q_{k}={kL}/
{2\pi }$ (see Figure \ref{fig:kse-solutions}(right)).  A truncation with $N \geq 8\nu =
8\floor{L/2\pi}$ can be considered accurate.  In the following, a
truncated system with $N=32\nu = 32\floor{L/2\pi}$ is considered to be
the ``full'' system, and we aim to construct reduced models
  for $K<2\nu \approx L/\pi$. Except when $N$ is small, the system
(\ref{DFM}) is stiff (since $q_{k}$ grows rapidly with $k$).  To handle
this stiffness and maintain reasonable accuracy, we generate data by
solving the truncated KSE by an exponential time difference fourth order
Runge-Kutta method (ETDRK4) \cite{CM02, KT05} with standard $3/2$
de-aliasing (see, e.g., \cite{Ors71, GO77}).  We solve the full system with a small step size $dt$, and then make observations of the modes with wave numbers $k=1,\dots ,K$ at a times separated by a
larger time interval $\delta > dt$, and denote the observed data by
\begin{equation*}
  \left\{ v_{k}(t_{n}),k=1,\dots ,K,n=0,\dots ,T\right\} ,
\end{equation*}
where $t_{n}=n\delta $.

To predict the evolution of $K$ observed modes with wave numbers
$k=1,\dots ,K$, it is natural to start from the truncated system that
includes only these modes, i.e., the system (\ref{DFM}) with $N=2(K+1)$.
However, when one takes $K$ to be relatively small (as we do in the present paper), large truncation errors are present, and the
dynamics of the truncated system are very different from those of the
full system.

\section{A discrete-time approach to stochastic parametrization and the
NARMAX representation}
\label{sect:narmax-review}

Suppose one is given a dynamical system
\begin{equation}
\frac{d\phi }{dt}=F(\phi )  \label{stoPara}
\end{equation}
where the variables are partitioned as $\phi =(u,w)$, with $u$
representing a (possibly quite small) subset of variables of direct
interest.  The problem of model reduction is to develop a reduced
dynamical system for predicting the evolution of $u$ alone.  That is,
one wishes to find an equation for $u$ that has the form
\begin{equation}
\frac{du}{dt}=R({u})+z(t),  \label{main}
\end{equation}
where $R(u)$ is a function of $u$ only, and $z(t)$ represents the model
error, i.e., the quantity one should add to $R(u)$ to obtain the correct
evolution. 
 For example, if $\phi$ represents the full state $(v_{k},k\in
\mathbb{Z})$ of a KSE solution and $u$ represents the low-wavenumber
modes $(v_{-K},v_{-K+1}, \cdots ,v_{K})$, then $R$ can correspond to a
$K$-mode truncation of the KSE.  In general, the model error $z$ must
depend on $u$, since the resolved variables $u$ and unresolved variables
$w$ typically interact.  

The usual approach to stochastic parametrization and model reduction as
formulated above is to identify $z$ as a stochastic process in the differential
equation (\ref{main}) from data (see \cite{Pal01, Wil05,CVE08,MH13} and
references therein). This approach has major difficulties. First, it leads to the
challenging problem of statistical inference for a continuous-time
nonlinear stochastic system from partial discrete observations
\cite{PSW09,ST12}.  The data are measurements of $u$, not of $z$; to
find values of $z$ one has to use equation (\ref{main}) and
differentiate $x$ numerically, which may be inaccurate because $z$ may
have high-frequency components or fail to be sufficiently smooth, and
because the data may not be available at sufficiently small time
intervals. Then, if one can successfully estimate values of $z$ and then
identify it, equation (\ref{main}) becomes a nonlinear stochastic
differential system, which may be hard to solve with sufficient accuracy
(see e.g \cite{KP99,MT04}).

To avoid these difficulties, a purely discrete-time approach to
stochastic parametrization was proposed in \cite{CL15}. This approach
avoids the difficult detour through a continuous-time stochastic system
followed by its discretization, by working entirely in a discrete-time
setting. It starts from the truncated equation 
\begin{equation*}
  \frac{dy}{dt}=R({y}).
\end{equation*}
($y$ differs from $u$ in that its evolution equation is missing the
information represented by the model error $z(t)$ in Eq.~(\ref{main}),
and so is not up to the task of computing $u$.)  Fix a step size
$\delta>0$, and choose a method of
time-discretization, for example fourth-order Runge-Kutta.  Then discretize
the truncated equation above to obtain a discrete-time approximation of
the form
\begin{equation*}
  y^{n+1}=y^{n}+\delta R^{\delta}(y^{n}).
\end{equation*}
(The function $R^\delta$ depends on the numerical time-stepping scheme used.)
To estimate the model error, write a discrete analog of equation
(\ref{main}):
\begin{equation}
u^{n+1}=u^{n}+\delta R^{\delta }(u^{n})+\delta z^{n+1}.  \label{uz}
\end{equation}
Note that a sequence of values of $z^{n+1}$ can be computed from data
using
\begin{equation}
z^{n+1}=\frac{u^{n+1}-u^{n}}{\delta }-R^{\delta }(u^{n}),  \label{zu}
\end{equation}
where the values of $u$ are the observed values; the resulting values of $z$
account for both the model error in (\ref{main}) and the numerical error in
the discretization $R^{\delta }(u^{n})$.  The task at hand is to identify the
time series $\left\{ z^{n}\right\} $ as a discrete stochastic process which depends
on $u$. Once this is done,
equation (\ref{uz}) will be used
to predict the evolution of $u$. There is no
need to approximate or differentiate, and there is no stochastic
differential equations to solve. Note that the $z^n$ depend on the numerical error as well as on the model error, and may not be good representations of the 
continuum model error; we are not interested in the latter, we are only interested
in modeling the solution $u$. 

The sequence $\left\{ z^{n}\right\} $ is a
stationary time series, which we represent via a NARMAX representation,
with $u$ as an exogenous input. This representation makes it possible to
take into account efficiently the non-Markovian features of the reduced
system as well as model and numerical errors.  The NARMAX representation
is versatile, easy to implement and reliable.  The model inferred from
data is exactly the same as the one used for prediction, which is not
the case for a continuous-time system because of numerical
approximations. The disadvantage of the discrete approach is that the
discrete system depends on both the spacing of the observed data and
  the method of time-discretization, so that 
data sets with different spacing lead to different discrete systems.

The NARMAX representation has the form:
\begin{eqnarray}
u^{n+1} &=&u^{n}+\delta R^{\delta }(u^{n})+\delta z^{n+1}, \notag \\
z^{n} &=&\Phi^{n}+\xi ^{n},  \label{z_equ}
\end{eqnarray}
for $n=1,2,\dots $, where $\left\{ \xi ^{n}\right\} $ is a sequence of
independent identically distributed random variables, the first equation repeats equation (\ref{uz}), and $\Phi ^{n}$
is a functional of current and past values of $\left( u,z,\xi \right) $,
of the parametrized form:
\begin{equation}
\Phi ^{n}=\mu +\sum_{j=1}^{p}A_{j}z^{n-j}+\sum_{j=1}^{r}B_{j}Q_{j}({u}^{n-j})+\sum_{j=1}^{q}C_{j}\xi ^{n-j},  \label{NARMAX}
\end{equation}
where $\left( Q_{j},\ i=1,\dots ,r\right) $ are functions to be chosen appropriately
and $\left( \mu ,A_{j},B_{j},C_{j}\right) $ are constant parameters to be
inferred from data. Here we assume that the real and complex parts of $\xi
^{n}$ are independent and have Gaussian distributions with mean zero and
diagonal covariance matrix $\sigma ^{2}$.

We call the above equations a NARMAX representation because the
time series $\{z^{n}\}$ in equations (\ref{z_equ})--(\ref{NARMAX})
resemble the nonlinear autoregression moving average with
exogenous input (NARMAX) model in e.g., \cite{Bil13,Han76,FY03}.
In (\ref{NARMAX}), the terms in $z$ are the autoregression part of order
$p$, the terms in $\xi $ are the moving average part of order $q$, and
the terms $Q_{j}$ which depend on $u$ are the exogenous input
terms. One should note that our representation is not quite a
NARMAX model in the usual sense, because the exogenous input in NARMAX is supposed to be independent of the output $z$, while here it
is not. The system (\ref{z_equ})--(\ref{NARMAX}) is a nonlinear
autoregression moving average (NARMA) representation of the time series
$\left\{ u^{n}\right\} $: by substituting (\ref{zu}) into (\ref{NARMAX})
we obtain a closed system for $u:$
\begin{equation}
  u^{n+1}=u^{n}+\delta R^{\delta }(u^{n})+\delta \Psi ^{n}+\delta \xi ^{n+1},
  \label{narma}
\end{equation}
where $\Psi ^{n}$ is a functional of the past values of $u$ and $\xi $.

The form of the system (\ref{z_equ})--(\ref{NARMAX}) is quite
general. $\Phi ^{n}$ can be a more general nonlinear functional of
the past values of $\left( u,z,\xi \right) $ than the one in
(\ref{NARMAX}). The additive noise $\xi $ can be replaced by a
multiplicative noise. We leave these further generalizations to
future work.  The main task in identification of a NARMAX representation is to determine the structure of the
functional $\Phi ^{n}$, that is, determine the terms that are needed,
then determine the orders $(p,r,q)$ in (\ref{NARMAX}) and estimate the parameters.

\section{Determination of the NARMAX representation}
\label{sect:narmax-struct}

To apply NARMAX to the KSE,
one must first decide on the structure of the model for $z^n$.
In particular, one must decide which nonlinear terms to include in the
{\em ansatz.}  We note that {\em a priori}, it is clear that {\em some}
nonlinear terms are necessary, as a simple linear {\em ansatz} is very
unlikely to be able to capture the dynamics of the KSE.  This is because
distinct modes are uncorrelated (see section 2), so that if the {\em
  ansatz} for $z^n$ contains only linear terms, then the different
components of our stochastic model for $z^n$ would be independent, and
one would not expect such a model to capture the energy balance between
the unresolved modes (see e.g. \cite{MH13}). 

But the question remains: which nonlinear terms? One option is to
include {\em all} nonlinear terms up to some order. This leads to a
difficult parameter estimation problems because of the large numbers of
parameters involved. Moreover, a model with many terms is likely to overfit, which complicates the model and leads to poor predictive performance (see e.g. \cite{Bil13}). What we do here instead is use inertial manifolds as a guide to selecting suitable nonlinear terms. We note that our construction satisfies the physical realizability constraints discussed in \cite{MH13}.

\subsection{Approximate inertial manifolds}

We begin with a quick review of approximate inertial
manifolds, following \cite{JKT90}.  Write the KSE (\ref{KSE}) in the form
\begin{equation*}
\frac{dv}{dt}=Av+f(v),\ \ v(x,0)=v_{0}\in H,
\end{equation*}
where the linear operator $A$ is the fourth derivative operator $-{
  \partial ^{4}}/{\partial x^{4}}$ or, in Fourier variables, the
infinite-dimensional diagonal matrix with entries $-q_{k}^{4}$ in
(\ref{FM}), and where $f(v)=-v{\partial v}/{\partial x}-{\partial
  ^{2}v}/{\partial x^{2}}$ or its spectral analog, and where $H=\left\{
v\in L^{2}[0,L]|v(x)=v(x+L),x\in \mathbb{R}\right\} $. An
  inertial manifold is a positively-invariant finite-dimensional
  manifold that attracts all trajectories exponentially;
%
%
see e.g. \cite{CFNT88, Rob94} for its existence, and see
e.g. \cite{TW94,JRT00} for estimates of its dimension. It can be
  shown that inertial manifolds for the KSE are realizable as graphs of
  functions $\psi:PH\rightarrow QH$, where $P$ is a suitable
  finite-dimensional projection (see below), and $Q=I-P$.  The inertial
  form (the ordinary differential equation that describes motion restricted to the
  inertial manifold) can be expressed in terms of the projected variable
  $u$ by 
\begin{equation}
\frac{du}{dt}=PAu+Pf(u+\psi (u)),  \label{pIF}
\end{equation}
where $u=Pv$.  Different methods for the approximation of inertial
manifolds have been proposed \cite{JKT90,JKT91, Ros95}. These methods
approximate the function $w=\psi (u)$ by approximating the projected
system on $QH$,

\begin{equation}
\frac{dw}{dt}=QAw+Qf(u+w),\ \ \ w(0)=Qv_{0}.  \label{w_equ}
\end{equation}
In practice, $P$ is typically taken to be the projection onto the
span of the first $m$ eigenfunctions of $A$, for example, one can set
$u=Pv=(v_{-m},v_{-m+1},\dots ,v_{m})$, the first $2m$ Fourier modes.
The dimension of the inertial manifold is generally not known {\em
    a priori} \cite{JRT00}, so $m$ is usually taken to be a large
integer.  It is shown in \cite{FMT88} that for large enough $m$ and for
$v=u+w$ on the inertial manifold, ${dw}/{dt}$ is relatively
small. Neglecting ${dw}/{dt}$ in (\ref{w_equ}) we obtain an
approximate equation for $w$:
\begin{equation}
w= -A^{-1}Qf(u+w).  \label{w_equ2}
\end{equation}
Approximations of $\psi $ can then be obtained by setting up 
fixed point iterations,
\begin{equation}
\psi _{0}=0,\,\,\,\psi _{n+1}=-A^{-1}Qf(u+\psi _{n}).  \label{fixed_pt}
\end{equation}
and stopping at a particular finite value of $n$, which yields an 
 ``approximate inertial manifold" (AIM). The accuracy of AIM improves as
its dimension $m$ increases, in the sense that the distance between the
AIM and the global attractor of the KSE decreases at the rate
$|q_{m}|^{-\gamma }$ for some $\gamma >0$ \cite{JKT90}.

In the application to our reduced KSE equation, one may try to set
$m=K$.
%
%
However, our cutoff $K<2\nu\approx \frac{L}{\pi }$ is too
  small for there to be a well-defined inertial form, since this is
  relatively close to the number $\nu=\floor{L/2\pi}$ of linearly
  unstable modes, and we generally expect any inertial manifold to have
  dimension much large than $2\nu$.  Nevertheless, we will see that the above
procedure for constructing AIMs provides a useful guide for selecting
nonlinear terms for NARMAX.  This is unsurprising, since inertial manifolds arise from the nonlinear energy transfer between low and high modes and the strong dissipation at high modes, which is exactly what we hope to capture.

\subsection{Structure selection}

We now determine the structure of the NARMAX model, i.e., decide what
terms should appear in (\ref{NARMAX}).  These terms should correctly
reflect how $z$ depends on $u$ and $\xi $. Once the terms are
determined, parameter estimation is relatively easy, as we see in the
next subsection.  The number of terms should be as small as possible, to
save computing effort and to reduce the statistical error in the
estimate of the coefficients.  One could try to determine the terms
to be kept by looking for the terms that contribute most to the output
variance, see \cite[Chapter 3]{Bil13} and the references there. This
approach fails to identify the right terms for strongly nonlinear
systems such as the one we have. Other ideas based on a purely
statistical approach have been explored e.g. in \cite[Chapter 9]{Bil13}.
In the present paper, we design the nonlinear terms in (\ref{NARMAX}) in the NARMAX model, which represents the model error in the $K$-mode truncation of the KSE,
using the theory of inertial manifolds sketched above.

In our setting, $u=\left( v_{-K},\dots ,v_{K}\right) $.
Because of the quadratic nonlinearity, only the modes with wave
number $|k|=K+1,K+2,\dots ,2K$ interact directly with the observed modes
in $u$; hence we set $w=\left( v_{-2K},\dots v_{-K-1},v_{K+1},\dots
v_{2K}\right) $ in Eq.~(\ref{w_equ2}). Using the result of a one-step
iteration $\psi _{1}=-A^{-1}Qf(u)$ in (\ref{fixed_pt}) we obtain an
expression for the high mode $v_k$ as a function of the low modes $u$:

\begin{equation*}
 \psi _{1,k}=- \left( A^{-1}Qf(u)\right) _{k}=- \frac{i}{2}q_{k}^{-4}\sum_{1\leq |l|\leq K,1\leq \left\vert
k-l\right\vert \leq K}v_{l}v_{k-l} ,
\end{equation*}
for $\left\vert k\right\vert =K+1,K+2,\dots ,2K$.

Our goal is to approximate the model error of the $K$-mode Galerkin truncation, $Pf(u+\psi (u))- Pf(u)$, so that the reduced model is closer to the attractor than the Galerkin truncation. In standard approximate inertial manifold methods, the model error is approximated by $Pf(u+\psi_1 (u))- Pf(u)$. Since $K$ is relatively small in our setting, we do not
   expect the AIM approximation to be effective.  However,
   in stochastic parameterization, we only need a parametric
   representation of the model error, and more specifically,
   to derive the nonlinear terms in the NARMAX model.  As
   the AIM procedure implicitly takes into account energy
   transfer between Fourier modes due to nonlinear
   interactions as well as the strong dissipation in high
   modes, it can be used to suggest nonlinear terms to
   include in our NARMAX model.  That is, the above
   expressions of $\psi_1$, $Pf(u+\psi_1 (u))- Pf(u)$
   provide explicit nonlinear terms $\{\phi_{1,k}
   \phi_{1,l}, \phi_{1,k} v_j\}$ we need (with
   $K<|k|,|l|\leq 2K, |j|\leq K$); we simply leave the
   coefficients of these terms as free parameters to be
   determined from data.  Roughly speaking, the NARMAX can
   be viewed as a regression implementation of a parametric
   version of AIM. 

Implementing the above observations in discrete time, we obtain the
following terms to be used in the discrete
reduced system:
\begin{equation}
\widetilde{u}_{j}^{n}=\left\{
\begin{array}{ll}
u_{j}^{n}, & 1\leq j\leq K; \\
i\sum_{l=j-K}^{K}u_{l}^{n}u_{j-l}^{n}, & K<j\leq 2K.
\end{array}
\right.  \label{u_ext}
\end{equation}
The modes with negative wave numbers are defined by $\widetilde{u}_{-j}^{n}=\widetilde{u}_{j}^{n,\ast }$. This yields the discrete
stochastic system
\begin{eqnarray*}
u_{k}^{n+1} &=&u_{k}^{n}+\delta R_{k}^{\delta }(u^{n})+\delta z_{k}^{n+1}, \\
z_{k}^{n} &=&\Phi _{k}^{n}+\xi _{k}^{n},
\end{eqnarray*}
where the
functional $\Phi _{k}^{n}$ has the form:
\begin{eqnarray}
\Phi _{k}^{n}(\theta _{k}) &=&\mu
_{k}+\sum_{j=1}^{p}a_{k,j}z_{k}^{n-j}+\sum_{j=0}^{r}b_{k,j}u_{k}^{n-j}+
\sum_{j=1}^{K}c_{k,j}\widetilde{u}_{j+K}^{n}\widetilde{u}_{j+K-k}^{n}  \notag
\\
&&+c_{k,(K+1)}R_{k}^{\delta }(u^{n})+\sum_{j=1}^{q}d_{k,j}\xi _{k}^{n-j},
\label{Phi}
\end{eqnarray}
for $1\leq k\leq K$. Here $\theta _{k}=\left( \mu
_{k},a_{k,j},b_{k,j},c_{k,j},d_{k,j}\right) $ are real parameters to be
estimated from data. Note that each one of the random variables $\xi_k^n$ affects directly only one mode $u_k^n$, as is consistent with the
  vanishing of correlations between distinct Fourier modes (see
  section~\ref{sect:kse}).

We set $\Phi _{-k}^{n}=\Phi _{k}^{n,\ast }$ so that the solution of the
stochastic reduced system satisfies $u_{-j}^{n}=u_{j}^{n,\ast }$.
We include the terms in $R^{\delta}_k(u^n)$ because the way they were introduced in the above reduced stochastic system does not guarantee that they have the optimal coefficients for the representation of $z$; the inclusion of these terms in $\Phi^n_k$ makes it possible to optimize these coefficients. This is similar in spirit to the construction of consistent reduced models in \cite{Har16}, though simpler to implement.


\subsection{Parameter estimation}
\label{section_paraEst}

We assume in this section that the terms and the orders $\left(
p,r,q\right) $ in the NARMAX representation have been selected, and
estimate the parameters 
as follows. To start, assume that the reduced system has dimension $K$,
that is, the variables $u^{n}$, $z^{n}$, $\Phi ^{n}$ and $\xi ^{n}$ have
$K$ components. Denote by $\theta _{k}=(\mu
_{k},A_{k,j},B_{k,j},C_{k,j})$ the set of parameters in the $k$th
component of $\Phi ^{n}$, and $\theta =(\theta _{1},\theta _{2},\dots
,\theta _{K})\,$. We write $\Phi _{k}^{n}$ as $\Phi _{k}^{n}(\theta_{k})$
to emphasize that $\Phi$ depends on $\theta _{k}$.


Recall that the real and complex parts of the components of $\xi ^{n}$
are independent $N(0,\sigma _{k}^{2})$ random variables. Then, following
(\ref{z_equ}), the log--likelihood of the observations $\left\{
u^{n},q+1\leq n\leq N\right\} $ conditioned on $\{\xi ^{1},\dots ,\xi
^{q}\}$ is (up to a constant)
\begin{equation}
l(\theta ,\sigma ^{2}|\xi ^{1},\dots ,\xi ^{q})=-\sum_{k=1}^{K}\left(
\sum_{n=q+1}^{N}\frac{\left\vert z_{k}^{n}-\Phi _{k}^{n}(\theta )\right\vert
^{2}}{2\sigma _{k}^{2}}+(N-q)\ln \sigma _{k}^{2}\right) .  \label{lkhd}
\end{equation}
If $q=0$, this is the standard likelihood of the data $\left\{u^{n},1\leq n\leq N\right\}$, and the values of $z^n$ and $\Phi^{n}(\theta)$ can be computed from the data $u^{n}$ using (\ref{zu}) and (\ref{Phi}), respectively. However, if $q>0$, the sequence $\{\Phi _{k}^{n}(\theta)\}$ cannot be computed directly from data, due to its dependence on the noise sequence $\{\xi^{n}\}$, which is unknown. Note that once the values of $\{\xi ^{1},\dots ,\xi^{q}\}$ are available, one can compute recursively the sequence $\{\Phi _{k}^{n}(\theta)\}$ for $n\geq q+1$ from data.  Hence we can compute the likelihood of  $\left\{u^{n},q+1\leq n\leq N\right\}$ conditional on $\{\xi ^{1},\dots ,\xi^{q}\}$.
If the stochastic reduced system is ergodic and the data come from this system, the MLE is asymptotically consistent (see e.g. \cite{Ham94, Han76}), and hence the values of $\xi ^{1},\dots ,\xi
^{q}$ do not affect the result if the data set is long enough. In practice, we can simply set $\xi ^{1}=\dots =\xi ^{q}=0$, the mean of these variables.

Taking partial derivatives with respect to $\sigma _{k}^{2}$ and noting
that $\left\vert z_{k}^{n}-\Phi _{k}^{n}(\theta _{k})\right\vert $ is
independent of $\sigma _{k}^{2}$, we find that the maximum likelihood
estimators (MLE) $\hat{\theta},\hat{\sigma}^{2}$ satisfy the following
equations:
\begin{equation}
\hat{\theta}_{k}=\arg
\min_{\theta_{k}}S_{k}(\theta_{k}),~\hat{\sigma}_{k}^{2}=\frac{1}{2(N-q)}S(\hat{\theta}_{k}),
\label{mle}
\end{equation}
where
\begin{equation*}
S_{k}(\theta _{k}):=\sum_{n=q+1}^{N}\left\vert z_{k}^{n}-\Phi
_{k}^{n}(\theta _{k})\right\vert ^{2}.
\end{equation*}

Note first that in the case $q=0$, the MLE $\hat{\theta}_{k}$ follows
directly from least squares, because $\Phi _{k}^{n}$ is linear in the
parameter $\theta _{k}$ and its terms can be computed from data.
If $q>0$,  the MLE $\hat{\theta}_{k}$ can be computed 
either by an optimization method  (e.g. quasi-Newton method), or by iterative
least squares method (see e.g. \cite{DC05}). With either method, one
first computes $\Phi _{k}^{n}(\theta _{k})$ with the current value of
parameter $\theta _{k}$, and then one updates $\theta _{k}$ (by gradient
search methods or by least squares), and repeats until the error
tolerance for convergence is reached. The starting values of $\theta _{k}$ for the iterations for
either method are set to be the least square estimates using the residual of the corresponding $q=0$ model. 

The simple forms of the log-likelihood in (\ref{lkhd}) and the MLEs in
(\ref{mle}) are based on the assumption that the real and complex parts
of the components of $\xi ^{n}$ are independent Gaussians. One may allow
the components of $\xi ^{n}$ to be correlated, at the cost of
introducing more parameters to be estimated. Also, similar to
\cite{DC05}, this algorithm can be implemented online, i.e. recursively as the data size increases,  and the noise sequence $\{\xi^n\}$ can be allowed to be non-Gaussian (on the basis
of martingale arguments).

\subsection{Order selection}

\label{section_order}

In analogy to the earlier discussion, it is not advantageous to have large orders
$(p,r,q)$ \cite{Bil13, BD02}, because, while large orders
generally yield small noise variance, the errors arising from the
estimation of the parameters accumulate as the number of parameters
increases. The forecasting ability of the reduced model depends not only
on the noise variance but also on the errors in parameter
estimation. For this reason, a penalty factor is often introduced
discourage the fitting of linear models with too many parameters. Many
criteria have been proposed for linear ARMA models (see
e.g. \cite{BD02}). However, due to the nonlinearity of NARMA and NARMAX models, these criteria do not work for them.

Here we propose a number of practical, qualitative criteria for
selecting orders by trial and error. We first select orders, estimate
the parameters for these orders, and then analyze how well the estimated
parameters and the resulting reduced system reach our goals. The
criteria are:

\begin{enumerate}
\item The variance of the model error should be small.

\item The stochastic reduced system should be stable, and its long-term
  statistical properties should be well-defined (i.e., the reduced
    system should have a stationary distribution), and should agree
  with the data. Especially, the autocorrelation functions
  (which are computed by time averaging) of the reduced system and of
  the data should be close.

 \item The estimated parameters should converge as the size of the data set
   increases.

\end{enumerate}

These criteria do not necessarily produce optimal solutions. As in most
statistics problems, one is aiming at an adequate rather than a perfect
solution.

\section{Numerical results}

\subsection{Estimated NARMAX coefficients}

We now determine the coefficients and and the orders in the functional
$\Phi$ of equation (\ref{NARMAX}) in the case $L=2\pi /\sqrt{0.085}$,
$K=5$; for this choice of $L$, the number of linearly unstable modes is
$\nu= \floor{1/\sqrt{0.085}}=3$.  This setting is the same as in Stinis
\cite{Sti04}, up to a change of variables in the solution. We obtain
data by solving the full Eq.~(\ref{DFM}) with $N=32\floor{L/2\pi}$ and
with time step $dt=0.001$, and make observations of the first $K$ modes
with wave number $k=1,\dots ,K$, with time spacing $\delta=0.1$. As
initial value we take $v_{0}(x)=(1+\sin x)\cos x$. Recall that we denote
by $\left\{ v(t_{n})\right\} _{n=1}^{T}$ the observations of the $K$
modes.  We choose the length of the data set to be large enough so that
the statistics can be computed by time averaging.  The means and
variances of the real parts of the observed Fourier modes settle down
after about $5\times 10^{4}$ time units. Hence we drop the first
  $10^4$ time units, and use observations of the next $5\times 10^{4}$
  time units as data for inferring a reduced stochastic system; with
  (estimated) integrated autocorrelation times of the Fourier modes
  ranging from $\approx 10$ to $\approx 35$ time units, this amount of
  data should be sufficient for estimating the coefficients.  (Because
$\delta=0.1$, the length of data is $T=5\times 10^{5}$.)

We consider different orders $(p,r,q)$: $p=0,1,2;r=1,2;q=0,1$. We first estimate the parameters by the conditional
likelihood method described in Section \ref{section_paraEst}. Then we
numerically test the stability of the stochastic reduced system by
generating a long trajectory of length $T$, starting from an arbitrary point in the
data (for example $v(t_{20000})$). We then drop the orders 
that lead to unstable systems, and select, among the remaining orders, the ones  with
smallest noise variances.

\begin{table}[tbp]
\caption{The noise variances in the NARMAX with different orders $(p,r,q)$.}
\label{tab:sigma}\centering
\begin{tabular}{l|l|lllll}
\hline
$(p,r,q)$ & scale & $\sigma _{1}^{2}$ & $\sigma _{2}^{2}$ & $\sigma _{3}^{2}$
& $\sigma _{4}^{2}$ & $\sigma _{5}^{2}$ \\ \hline
010 & $\times 10^{-3}$ & 0.0005 & 0.0061 & 0.0217 & 0.1293 & 0.1638 \\
020 & $\times 10^{-3}$ & 0.0005 & 0.0051 & 0.0187 & 0.0968 & 0.1612 \\
110 & $\times 10^{-4}$ & 0.0047 & 0.0587 & 0.1664 & 0.4290 & 0.5640 \\
120 & $\times 10^{-4}$ & 0.0044 & 0.0509 & 0.1647 & 0.4148 & 0.3259 \\
021 & $\times 10^{-4}$ & 0.0012 & 0.0129 & 0.0472 & 0.2434 & 0.4056 \\
111 & $\times 10^{-4}$ & 0.0012 & 0.0148 & 0.0421 & 0.1079 & 0.1426 \\
210 & $\times 10^{-5}$ & 0.0020 & 0.0307 & 0.1533 & 0.3921 & 0.3234 \\
220 & $\times 10^{-5}$ & 0.0019 & 0.0304 & 0.1450 & 0.2858 & 0.1419 \\ \hline
\end{tabular}
\end{table}
The orders $(2,1,0),(2,2,0)$ lead to unstable reduced system, though
their noise variances are the smallest, see Table~\ref{tab:sigma}. This suggests that large orders $p,r,q$ are not needed.
Among the other orders, the choices $(0,2,1)$ and $(1,1,1)$ have
the smallest noise variances (see Table \ref{tab:sigma}). The orders $(1,1,1)$ seems to be better than $(0,2,1)$, because the
former has four out of the five variances smaller than the latter.

\begin{table}[tbp]
\caption{The average of mean square distances between the autocorrelation functions of the
  data and the autocorrelation functions NARMAX trajectory with
  different orders $(p,r,q)$. }
\label{tab:distance}\centering
\begin{tabular}{c|ccccc}
\hline
$(p,r,q)$ & $D_{1}$ & $D_{2}$ & $D_3$ & $D_{4}$ & $D_{5}$ \\ \hline
020 & 0.0012 & 0.0013 & 0.0008 & 0.0004 & 0.0009 \\
021 & 0.0010 & 0.0010 & 0.0007 & 0.0004 & 0.0008 \\
111 & 0.0013 & 0.0019 & 0.0010 & 0.0004 & 0.0010 \\
110 & 0.0011 & 0.0018 & 0.0008 & 0.0005 & 0.0010 \\
120 & 0.0013 & 0.0022 & 0.0010 & 0.0004 & 0.0011 \\ \hline
\end{tabular}
\end{table}
For further selection, following the second criterion in section
\ref{section_order}, we compare the empirical autocorrelation functions
of the NARMAX\ reduced system with the autocorrelation functions of
data. Specifically, first we take $N_{0}$ pieces of the data, $\left\{
\left( v(t_{n}),n=n_{i},n_{i}+1,\dots ,n_{i}+T\right) \right\}
_{i=1}^{N_{0}}$ with $n_{i+1}=n_{i}+T_{lag}/\delta$, where $T$ is the
length of each piece and $T_{lag}$ is the time gap between two adjacent
pieces. For each piece $\left( v(t_{n}),n=n_{i},\dots ,n_{i}+T\right) $,
we generate a sample trajectory of length $T$ from the NARMAX reduced
system using initial $\left( v(t_{n_{i}}),v(t_{n_{i}+1}),\dots
,v(t_{n_{i}+m})\right) $, where $m=\max \left\{ p,r,2q\right\} +1$, and
denote the sample trajectory by $\left( u^{n_{i}+n},n=1,\dots ,T\right)
$.  Here an initial segment is used
to estimate the first few steps of the noise sequence, $\left( \xi
^{q+1},\dots ,\xi ^{2q}\right) $ (recall that we set $\xi ^{1}=\dots
=\xi ^{q}=0\,$). Then we compute the autocorrelation functions of real
parts each trajectory by

\begin{eqnarray*}
\gamma _{v,k}(h,i) &=&\frac{1}{T-h}\sum_{n=1}^{T-h}\operatorname{Re}v_{k}(t_{n_{i}+n+h})\operatorname{Re}v_{k}(t_{n_{i}+n}); \\
\gamma _{u,k}(h,i) &=&\frac{1}{T-h}\sum_{n=1}^{T-h}\operatorname{Re}u_{k}^{n_{i}+n+h}\operatorname{Re}u_{k}^{n_{i}+n}
\end{eqnarray*}
for $h=1,\dots ,T_{lag}/\delta$, $i=1,\dots ,N_{0}$, $k=1,\dots ,K$,
and compute the average of mean square distances between the
autocorrelation functions by
\begin{equation*}
D_{k}=\frac{1}{N_{0}}\sum_{i=1}^{N_{0}}\left( \frac{\delta}{T_{lag}}
\sum_{h=1}^{T_{lag}/\delta}\left\vert \gamma _{v,k}(h,i)-\gamma
_{u,k}(h,i)\right\vert ^{2}\right) .
\end{equation*}
The orders $(p,r,q)$ with the smallest average mean square
distances will be selected. Here we only consider the autocorrelation functions of the real
parts, since the imaginary parts have statistical properties similar to those of the real
parts. 


Table \ref{tab:distance} shows the average of mean square
    distances between the autocorrelation functions for different orders and the autocorrelation functions of the data, computed with $N_{0}=100,T_{lag}=50$. The orders $(1,1,1)$ have
larger distances than $(0,2,1)$. We select the orders $(0,2,1)$
because they have the smallest average of mean square distances. The estimated
parameters for the orders $(0,2,1)$ are presented in Table
\ref{tab:021}.

\begin{table}[tbp]
\caption{The estimated parameters in the NARMAX with orders $(p,r,q)=(0,2,1)$.}
\label{tab:021}\centering
\begin{tabular}{c|rrrrrrr}
\hline
$k$ & $\mu _{k}$($\times 10^{-4}$) & $b_{k,0}$ & $b_{k,1}$ & $d_{k,1}$ & $\sigma _{k}^{2}$($\times 10^{-4}$) &  &  \\ \hline
1 & 0.0425 & 0.0909 & -0.0910 & 0.9959 & 0.0012 &  &  \\
2 & -0.0073 & 0.1593  & -0.1600 & 0.9962 & 0.0138 &  &  \\
3 & 0.3969 & 0.2598 & -0.2617  & 0.9942 & 0.0520 &  &  \\
4 & -0.9689 & 0.7374 & -0.7408  & 0.9977 & 0.2544 &  &  \\
5 & -0.1674& 0.3822 & -0.3799 & 0.9974 & 0.4056 &  &  \\ \hline\hline
$k$ & $c_{k,1}$ & $c_{k,2}$($\times 10^{-3}$) & $c_{k,3}$($\times 10^{-3}$)
& $c_{k,4}$($\times 10^{-3}$) & $c_{k,5}$ ($\times 10^{-3}$)& $c_{k,6}$  \\ \hline
1 & 0.0002 & 0.0000 & 0.0010 & 0.0013 & -0.0003 & -0.0082  \\
2 & 0.0005 & 0.2089 & 0.0015 & 0.0007 & 0.0001 & -0.0157  \\
3 & 0.0008 & 0.3836 & 0.1055 & -0.0013 & -0.0040 & -0.0283  \\
4 & 0.0010 & 0.5841 & 0.2971 & -0.4104 & 0.0449 & -0.0800  \\
5 & 0.0012 & 0.6674 & 0.4763 & 0.2707 & 0.1016 & -0.0710 \\ \hline
\end{tabular}
\end{table}

For comparison, we also carried out a similar analysis for
$\delta=0.01$.  We found (data not shown) that (i) the best choices of
$(p,r,q)$ are different for $\delta=0.1$ and for $\delta=0.01$;
and (ii) the coefficients do not scale in a simple way, e.g., all as
some power of $\delta$.  Presumably, there is an asymptotic regime (as
$\delta\to0$) in which the coefficients do exhibit some scaling, but
$\delta=0.1$ is too large to exhibit any readily-discernible scaling
behavior.  We leave the investigation of the scaling limit as
$\delta\to0$ for future work.

\medskip

We observed that for all the above orders, most of the
estimated parameters show a clear trend of convergence as the length of the
data set increases, but some parameters keep oscillating (data not shown). For example, for NARMAX with orders $(0,2,1)$, the coefficients $b_{k,j}$ are more oscillatory than
the coefficients $c_{k,j}$, and the parameters of the mode with wave number $k=5$ are more oscillatory than the parameters of other modes. This indicates
that the structure and the orders are not yet optimal, and we
leave the task of developing better structures and orders to future work.
Here we select the orders simply based on the size of noise variances and
on the ability to reproduce the autocorrelation functions. Yet we obtain a reduced system which
achieves both our goals of reproducing the long-term statistics and making
reliable short-term forecasting, as we show in the following sections.

In the following, we select the orders $(0,2,1)$ for the NARMAX
reduced system.

\subsection{Long-term statistics}
We compare the statistics of the truncated system and the NARMAX reduced system with the statistics
of the data. We calculate the following quantities for the reduced systems as
well as for data: the empirical probability density functions (pdf) and the empirical autocorrelation functions for each of the $K$
components. All these statistics are computed by time-averaging long sample
trajectories, as we did for the autocorrelation functions in the
previous subsection.

\begin{figure}
  \centering
    \includegraphics[width=0.95\textwidth]{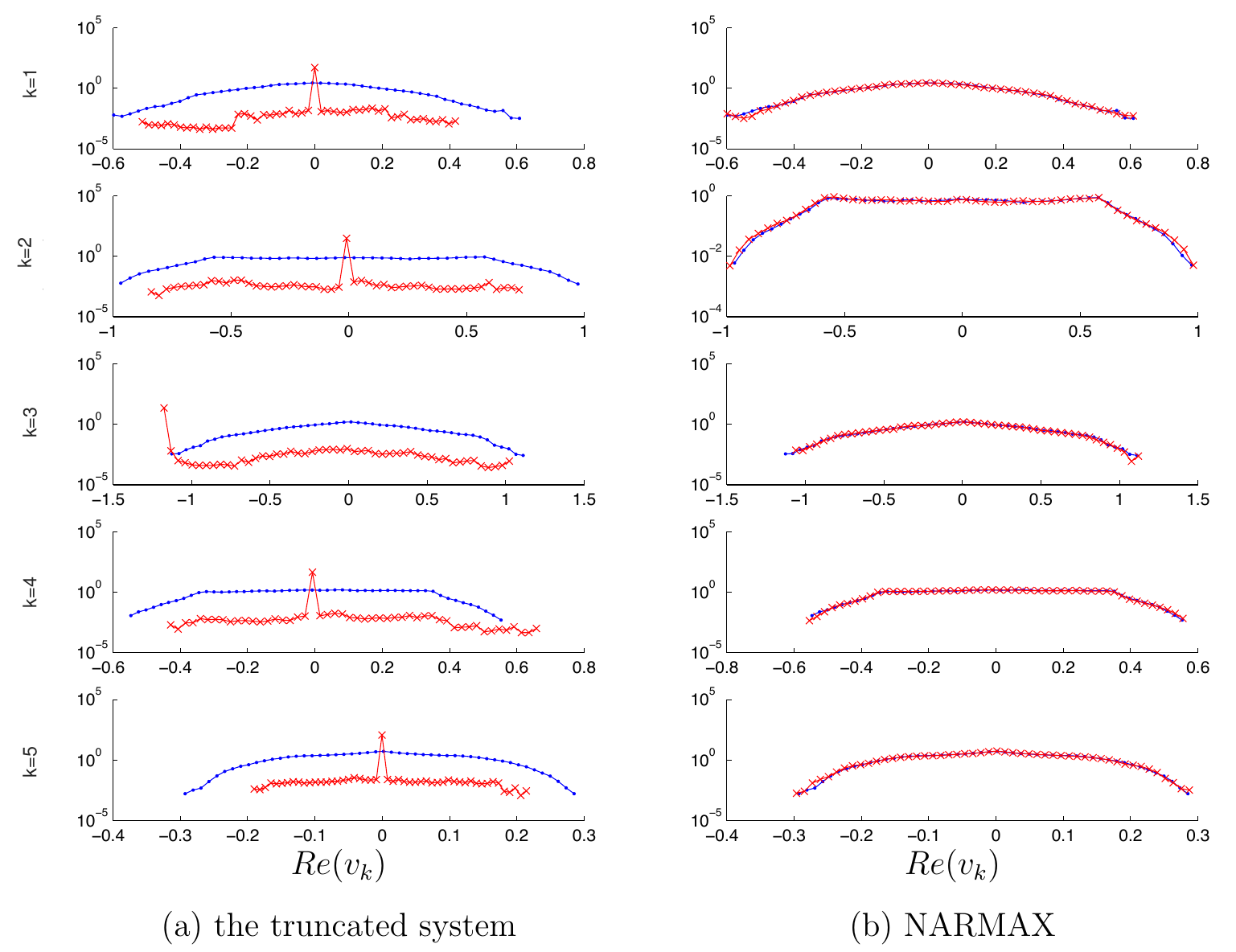} 
    \caption{Semilog plot of the probability density functions
      reproduced by the truncated system and by the NARMAX system (red
      $\times$'s), compared to the full system (blue dots). Here we only
      plot the real parts of the paths of the Fourier modes (from top to
      bottom, the wave numbers are $k=1,2,\dots,5$).  Note the
      distributions of the truncated system are not symmetric because
      solutions appear to converge to a fixed point, so that our
      ergodicity assumption is invalid.}
    \label{fig:pdf}
\end{figure}

\begin{figure}
  \centering
\includegraphics[width=0.95\textwidth]{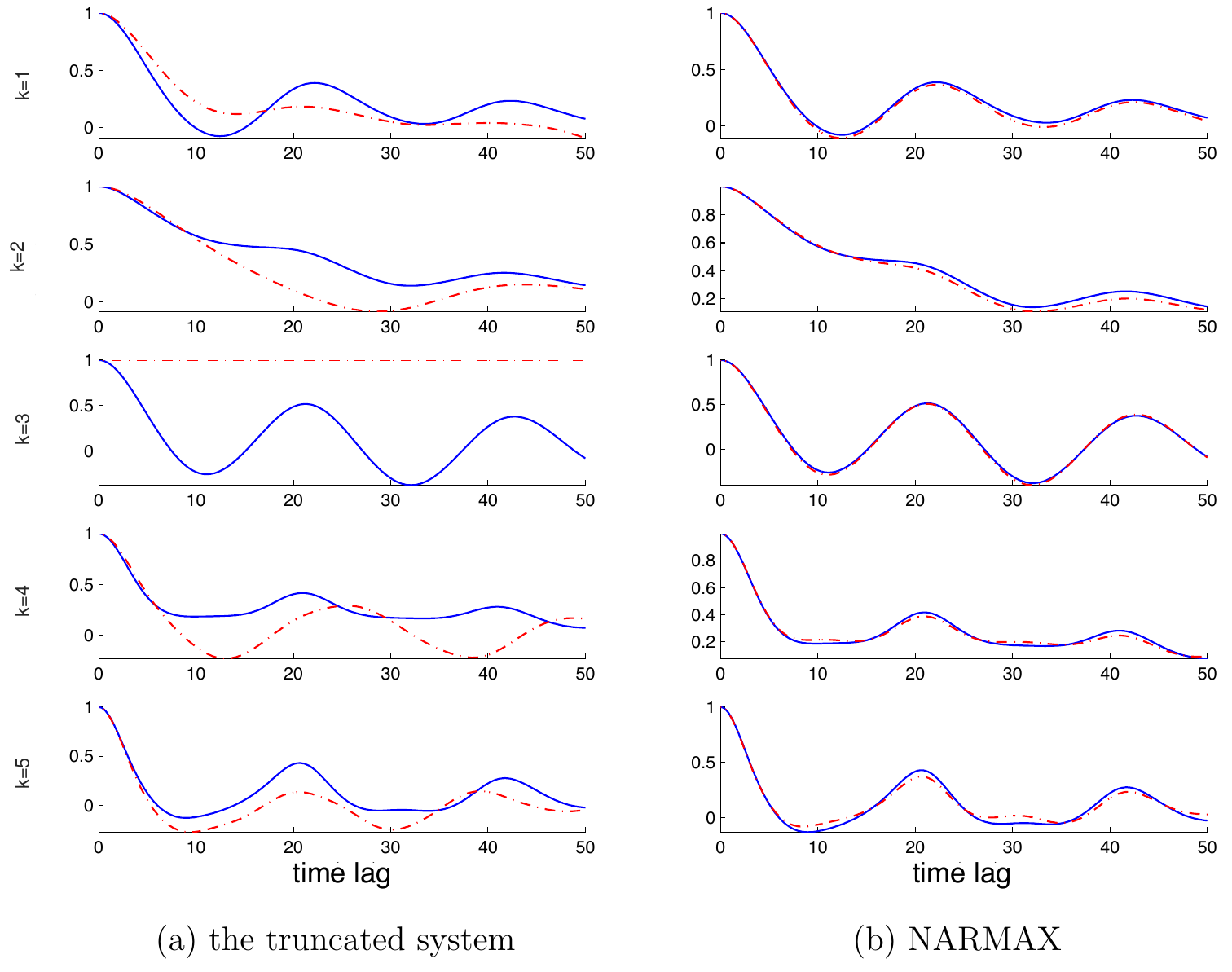} 
     \caption{The autocorrelation functions reproduced by the truncated
       system and and the NARMAX system (dot dash line), compared to the
       full system (solid line). Here we only plot the real parts of the
       paths of the Fourier modes (from top to bottom, the wave numbers
       are $k=1,2,\dots,5$).}
      \label{fig:acf}
 \end{figure}

The pdfs and autocorrelation functions of data are reproduced well by the NARMAX reduced
system, as shown in Figure \ref{fig:pdf} and Figure \ref{fig:acf}. The
NARMAX system reproduces almost exactly the pdfs and the autocorrelation functions of the data, a
significant improvement over the truncated system.

We also computed energy covariances ${\rm
    cov}(|v_k|^2,|v_\ell^2|)$ for some pairs $(k,\ell)$ with
  $k\neq\ell$.  (Recall from section~\ref{sect:kse} that this is
  expected to be nonzero.)  In particular, tests with $k=2$ and
  $\ell\in\{1,2,3,4,5\}$ show that the NARMAX reduced
  system can correctly capture these energy-energy covariances, whereas
  the truncated system does not (data not shown).

\subsection{Short-term forecasting}

We now investigate how well the NARMAX reduced system predicts the behavior
of the full system.

\begin{figure}
    \centering
        \includegraphics[width=1.0\textwidth]{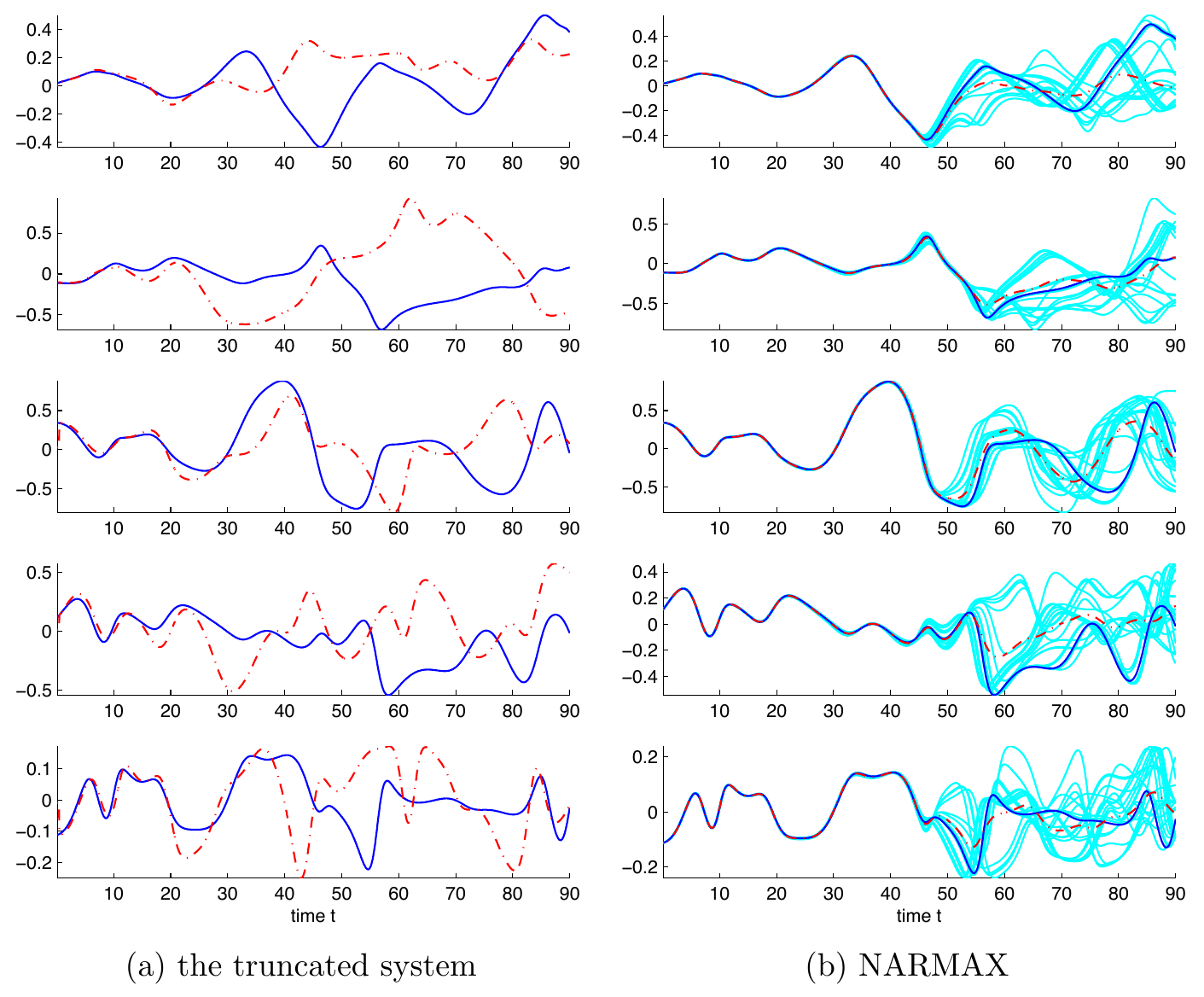}
    \caption{A typical path forecast made by the reduced systems
      compared to the data (solid line).  (a) the truncated system
      (dot-dashed line); (b) ensemble trajectories (cyan line) and the
      mean trajectory of the ensemble (dot-dashed line) of the NARMAX
      system.  Here we only plot the real parts of the paths of the
      Fourier modes (from top to bottom, the wave numbers are
      $k=1,\cdots,5$).} \label{fig:path0}
\end{figure}

We start from single path forecasts. For the NARMAX with orders
$(p,r,q)$, we start the multistep recursion by using an initial segment with $m=2\max\{p,r,q\}+1$ steps as follows. We set $\xi^1=\cdots=\xi^q=0$, and estimate
$\xi^{q+1},\dots,\xi^{m}$ using equation (\ref{z_equ}). Then we follow
the discrete system to generate an ensemble of trajectories from
different realizations, with all realizations using the same initial condition. We do not
introduce artificial perturbations into the initial conditions, because
the exact initial conditions are known. A typical ensemble of $20$
trajectories, as well as its mean trajectory and the corresponding data
trajectory, is shown in Figure \ref{fig:path0}(b). As a comparison, we
also plot a forecast using the truncated system. Since the observations
provide the exact initial condition, the truncated system produces a
single forecast path, see Figure \ref{fig:path0}(a). We observe that the
ensemble of NARMAX follows the true trajectory well for about 50 time
units, and the spread becomes wide quickly afterwards, while the
ensemble mean can follow the true trajectory to 55 time units. Compared
to the truncated system which can make forecast for about 20 time units,
NARMAX can make forecast for about 55 time units, which is a significant
improvement. We comment that in the prediction time is about 35 time
units in \cite{Sti04}, where the Mori-Zwanzig formalism is used.

To measure the reliability of the forecast as a function of lead time, we compute two commonly-used statistics, the root-mean-square-error (RMSE) and the anomaly correlation (ANCR).  Both statistics are based on generating a large number of ensembles of trajectories of the reduced system starting from different initial conditions, and comparing the mean ensemble predictions with the true trajectories, as follows.

First we take $N_{0}$ short pieces of the data, $\left\{ \left(
v(t_{n}),n=n_{i},n_{i}+1,\dots ,n_{i}+T\right) \right\} _{i=1}^{N_{0}}$ with
$n_{i+1}=n_{i}+T_{lag}/\delta$, where $T=$ $T_{lag}/\delta$ is the
length of each piece and $T_{lag}$ is the time gap between two adjacent
pieces. For each short piece of data $\left( v(t_{n}),n=n_{i},\dots
,n_{i}+T\right) $, we generate $N_{ens}$ trajectories of length $T$ from the
NARMAX reduced system, starting all ensemble members from the same
several-step initial condition $\left( v(t_{n_{i}}),v(t_{n_{i}+1}),\dots
,v(t_{n_{i}+m})\right) $, where $m=2\max \left\{ p,r,q\right\} +1$, and
denote the sample trajectories by $\left( u^{n}(i,j),n=1,\dots ,T\right) $
for $i=1,\dots ,N_{0}$ and $j=1,\dots ,N_{ens}$. Again, we do not introduce
artificial perturbations into the initial conditions, because the exact
initial conditions are known, and by initializing from data, we preserve
the memory of the system so as to generate better ensemble trajectories.

We then calculate the mean trajectory for each ensemble, $\bar{u}^{n}(i)=
\frac{1}{N_{ens}}\sum_{j=1}^{N_{ens}}u^{n}(i,j)$. The RMSE measures,  in an average sense, the
difference between the mean ensemble trajectory, i.e., the expected path predicted by the reduced model, and the true data trajectory:

\begin{equation*}
\mathrm{RMSE}(\tau _{n}):=\left( \frac{1}{N_{0}}\sum_{i=1}^{N_{0}}\left\vert \operatorname{Re}
v(t_{n_{i}+n})-\operatorname{Re}\bar{u}^{n}(i)\right\vert ^{2}\right)^{1/2},
\end{equation*} 
where $\tau _{n}=n\delta$. The anomaly correlation (ANCR) shows the
average correlation between the mean ensemble trajectory and the true data
trajectory (see e.g \cite{CVE08}):
\begin{equation*}
\mathrm{ANCR}(\tau _{n}):=\frac{1}{N_{0}}\sum_{i=1}^{N_{0}}\frac{\mathbf{a}
^{v,i}(n)\cdot \mathbf{a}^{u,i}(n)}{\sqrt{|\mathbf{a}^{v,i}(n)|^{2}\left
\vert \mathbf{a}^{u,i}(n)\right\vert ^{2}}},
\end{equation*}
where $\mathbf{a}^{v,i}(n)=\operatorname{Re}v(t_{n_{i}+n})-\operatorname{Re}\left\langle
v\right\rangle $ and $\mathbf{a}^{u,i}(n)=\operatorname{Re}\bar{u}^{n}(i)-\operatorname{Re}
\left\langle v\right\rangle $ are the anomalies in data and the ensemble
mean. Here $\mathbf{a\cdot b=}\sum_{k=1}^{K}a_{k}b_{k}$, $\left\vert \mathbf{
a}\right\vert ^{2}=\mathbf{a\cdot a}$, and $\left\langle v\right\rangle $ is
the time average of the long trajectory of $v$. Both statistics measure the accuracy of the mean ensemble prediction; $\rm{RMSE}=0$ and $\rm{ANCR}=1$ would correspond to a perfect prediction, and small RMSEs and large (close to 1) ANCRs are desired.

\begin{figure}[tbp]
\begin{center}
\begin{tabular}{cc}
 \resizebox{0.45\textwidth}{!}{\includegraphics{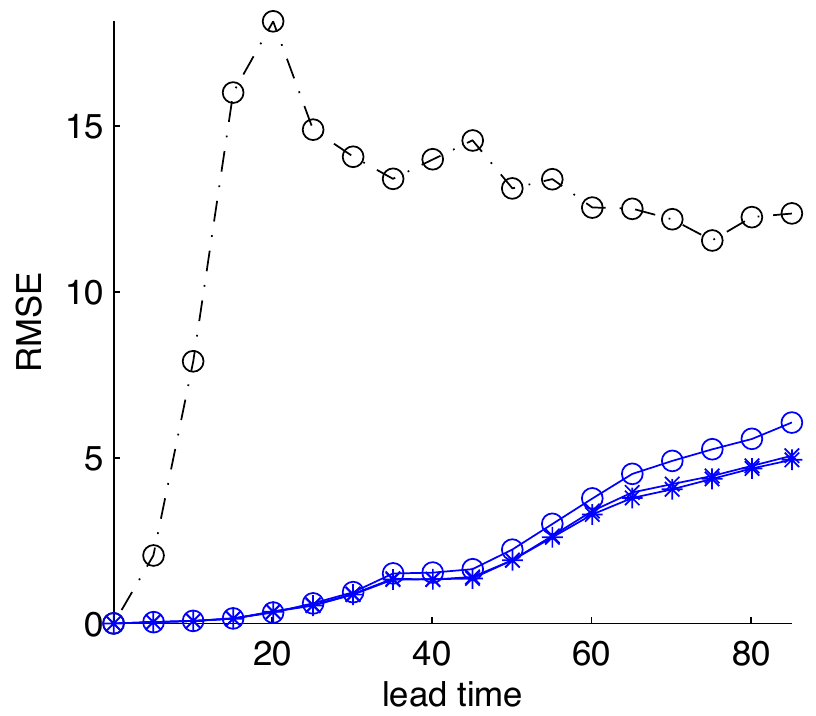}}  &
 \resizebox{0.45\textwidth}{!}{\includegraphics{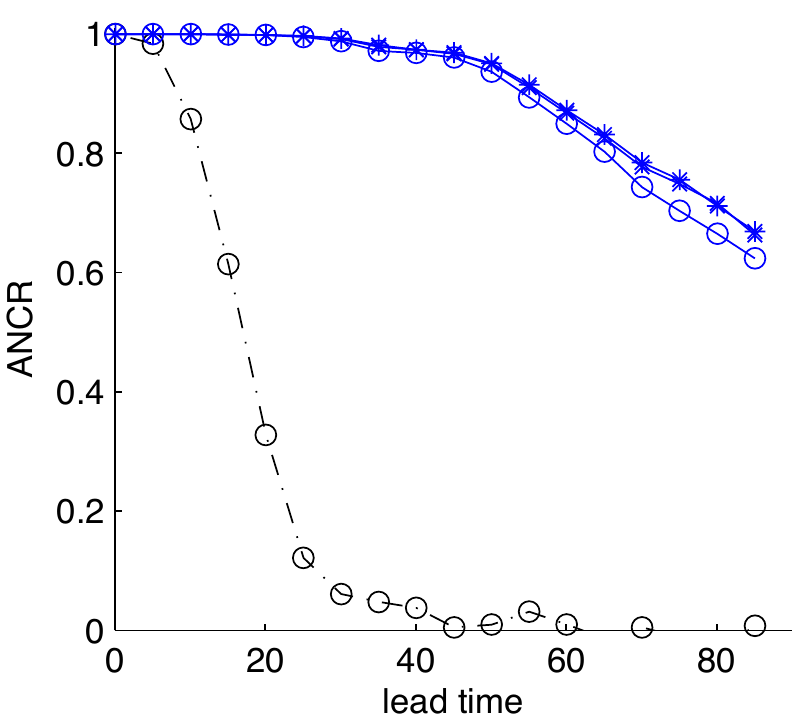}} \\ 
 (a) RMSE ~ &  (b) ANCR 
 \end{tabular}
\end{center}
\vspace{-4mm}
\caption{Root-mean-square-error (RMSE) and Anomaly correlations (ANCR) of ensemble forecasting, produced by the NARMAX system (solid lines) and the truncated system (dot dashed line), for different ensemble sizes: $N_{ens}=1$ (circle marker), $N_{ens}=5$ (cross marker), and $N_{ens}=20$ (asterisk marker).} \label{fig:RMSE}
\end{figure}

Results for RMSE and ANCR for $N_{0}=1000$ ensembles are shown in Figure \ref{fig:RMSE}, where we tested three ensemble sizes: $N_{ens}=1,5,20$. The forecast lead time at which the RMSE keeps below 2 is about 50 time units, which is about 10 times of the forecast lead time of the truncated system. The forecast lead time at which the ANCR drops below 0.9 is about 55 time units, which is about five times of number of the truncated system. We also observe that a larger ensemble size leads to smaller RMSEs and larger ANCRs.
\subsection{Importance of the nonlinear terms in NARMAX}

\begin{figure}
         \centering
        \includegraphics[width=0.9\textwidth]{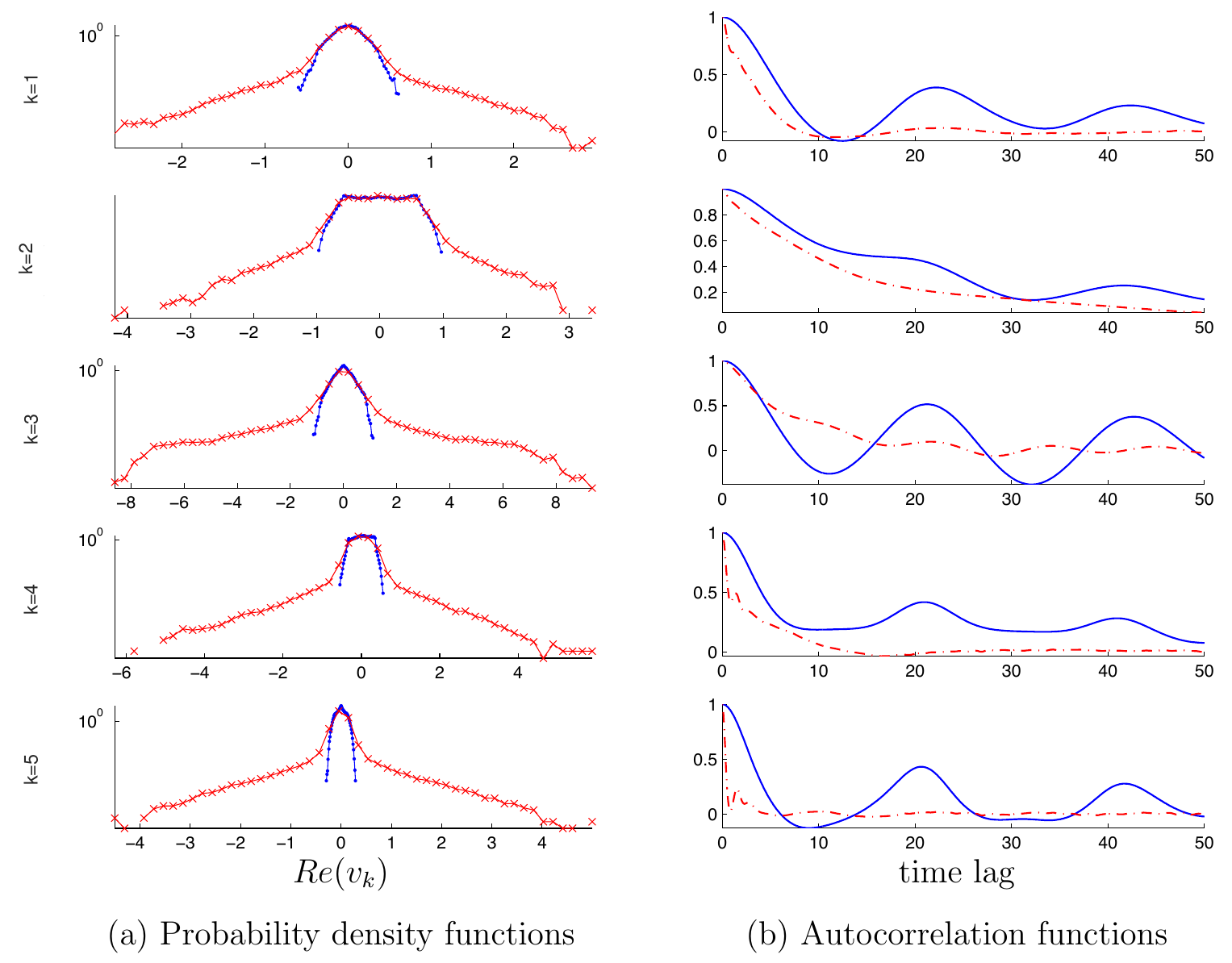}
  \caption{Results produced by ARMAX, i.e., using only linear terms in
    the ansatz for $z^n$.  In (a), we compare the probability
    density functions of different modes produced by ARMAX (red
    $\times$'s) to that of the full model (blue dots).  In (b), we
    compare the autocorrelation functions (solid line for full model,
    dashed line for ARMAX).  Compare to Figures~\ref{fig:pdf}(b) and
    \ref{fig:acf}(b).}
  \label{fig:armax}
\end{figure}

Finally, we examine the necessity of including nonlinear terms in the ansatz, by comparing NARMAX to ARMAX, i.e.,
stochastic parametrization keeping only the linear terms in the ansatz.  We performed a number of numerical experiments in which
we fitted an ARMAX ansatz to data.  We found that for $(p,r,q) =
(0,2,1)$, which was the best order we found for NARMAX, the
corresponding ARMAX approximation was unstable.

We also found the best orders for ARMAX to be stable, which were
$(p,r,q)=(2,1,0)$, and compared the results to those produced by NARMAX.
The results are shown in Figure~\ref{fig:armax}.  In (a), the pdfs are
shown for each resolved Fourier mode, and compared to those of the full
model.  Clearly, the Fourier modes for ARMAX experience much larger
fluctuations, presumably because of the build-up of energy in the
resolved modes.  In contrast, the results produced by NARMAX match
reality much better (see Figures~\ref{fig:pdf}(b)), as it more correctly
models the nonlinear interactions between different modes.  A
consequence of these larger fluctuations is that ARMAX cannot even
capture the mean energy spectrum: the ARMAX model leads to mean energies
that are about 5 times larger than the true energy spectrum, which
NARMAX is able to reproduce (data not shown).

Figure~\ref{fig:armax}(b) shows the corresponding autocorrelations.  We
see that ARMAX does not correctly capture the temporal structure of the
dynamics.  Again, this is consistent with the fact that in ARMAX the components of $z^n$ are independent, which do not correctly
capture the energetics the KSE.  In contrast, the NARMAX results in
Figure~\ref{fig:acf}(b) show a much better match.

\section{Conclusions and discussion}

We performed a stochastic parametrization for the KSE equation in order
to use it as a test bed for developing such parametrization for more
complicated systems.  We estimated and identified the model error in a
discrete-time setting, which made the inference from data easier, and
avoided the need to solve nonlinear stochastic differential systems; we then
represented the model error as a NARMAX time series.  We found an
efficient form for the NARMAX series with the help of an approximate
inertial manifold, which we determined by a construction developed in a
continuum setting, and which we improved by parametrizing its
coefficients.

A number of dimensional reduction techniques have been developed over
the years in the continuum setting, e.g., inertial manifolds,
renormalization groups, the Mori-Zwanzig formalism, and a variety of
perturbation-based methods. In the present paper we showed, in the
Kuramoto-Sivashinsky case, that continuum methods could be adapted for
use in the more practical discrete-time setting, where they could help
to find an effective structure for the NARMAX series, and could in turn
be enhanced by estimating the coefficients that appear, producing an
effective and relatively simple parametrization. Another example in a
similar spirit was provided by Stinis \cite{Sti15}, who renormalized
coefficients in a series implementation of the Mori-Zwanzig formalism.

Such continuous/discrete, analytical/numerical hybrids raise interesting
questions. Do there exist general, systematic ways to use continuum
models to identify terms in NARMAX series? Does the discrete setting
require in general that the continuum methods be modified or
discretized? What are the limitations of this approach? The answers
await further work.
\section{Acknowledgements}
The authors thank the anonymous referee, Prof. Panos Stinis, and Prof. Robert Miller for their helpful suggestions.
KL is supported in part by the National Science Foundation under grants DMS-1217065 and DMS-1418775, and thanks the Mathematics Group at Lawrence Berkeley National Laboratory for facilitating this
  collaboration.  AJC and FL are supported in part by the Director,
Office of Science, Computational and Technology Research,
U.S. Department of Energy, under Contract No. DE-AC02-05CH11231, and by
the National Science Foundation under grants DMS-1217065 and
DMS-1419044.

\bibliographystyle{plain}
\bibliography{references}

\end{document}